\documentclass[11pt]{article}
\usepackage{amssymb}
\newtheorem{Lemma}{Lemma}
\newtheorem{Theorem}{Theorem}
\newtheorem{Proposition}[Lemma]{Proposition}
\newtheorem{Remark}{Remark}

\newenvironment{Proof}%
 {\begin{trivlist} \item[]{\bf Proof. }}%
 {\hspace*{\fill}$\rule{.3\baselineskip}{.35\baselineskip}$\end{trivlist}}
 \newenvironment{Proof1}%
{\begin{trivlist} \item[]{\bf Proof }}%
{\hspace*{\fill}$\rule{.3\baselineskip}{.35\baselineskip}$\end{trivlist}}
 {\begin{trivlist}\item[]\textbf{Acknowledgments }}{\end{trivlist}}

\makeatletter \@addtoreset{equation}{section} \makeatother

\newcommand{\C}{\mathbb{C}}

\newcommand{\R}{\mathbb{R}}

\newfam\bifam
\font\tenbi=cmmib10 scaled \magstep1 \font\sevenbi=cmmib10 at 11pt
\font\fivebi=cmmib10 at 6pt \textfont\bifam = \tenbi
\scriptfont\bifam= \sevenbi \scriptscriptfont\bifam= \fivebi

\begin{document}

\title{\bf Spectral stability of nonlinear waves \\
in KdV-type evolution equations}

\author{Dmitry Pelinovsky \\
{\it Department of Mathematics, McMaster University, Hamilton, Ontario, Canada, L8S 4K1}}

\date{\today}
\maketitle

\begin{abstract}
This paper concerns spectral stability of nonlinear waves in
KdV-type evolution equations. The relevant eigenvalue problem
is defined by the composition of an unbounded self-adjoint operator with
a finite number of negative eigenvalues and an unbounded
non-invertible operator $\partial_x$. The instability index theorem is
proven under a generic assumption on the self-adjoint operator
both in the case of solitary waves and periodic waves.
This result is reviewed in the context of recent results
on spectral stability of nonlinear waves in KdV-type evolution equations.
\end{abstract}


\section{Introduction}

KdV-type evolution equations are defined by the following nonlinear PDE in $(1+1)$ variables:
\begin{equation}
\label{Ham-system} \frac{\partial u}{\partial t} = \frac{\partial}{\partial x} E'(u), \quad u(t) \in
{\cal X},
\end{equation}
where $E: \mathcal{X} \rightarrow \mathbb{R}$ is a $C^2$ functional on a subspace $\mathcal{X}$
of Hilbert space $L^2$ associated with the inner product $\langle \cdot,\cdot \rangle$ and an induced norm
$\| \cdot \|$.
A critical point $\phi \in {\cal X}$ of the Hamiltonian
functional $E$, defined by $E'(\phi) = 0$, represents a nonlinear wave of the KdV-type evolution
equation. Depending on the phase space $\mathcal{X}$, $\phi$ can be a solitary wave on an infinite
line $\mathbb{R}$ or a $(2L)$-periodic wave on the fundamental period $[-L,L]$.
In what follows before the last section, we consider solitary waves
on an infinite line defined in a $L^2$-based Sobolev space.

The spectral stability of $\phi$ is determined by the
spectrum of the non-self-adjoint eigenvalue problem
\begin{equation}
\label{spectrum} \partial_x E''(\phi) v = \lambda v,
\end{equation}
where $\mathcal{L} := E''(\phi)$ is a self-adjoint real-valued operator with
a dense domain $D(\mathcal{L})$ in $L^2(\mathbb{R})$. Since $\mathcal{L}$ is real-valued,
for every eigenvalue $\lambda \in \C$ with ${\rm Im}(\lambda) \neq 0$,
there is an eigenvalue $\bar{\lambda} \in \C$. We also assume the Hamiltonian symmetry,
that is, for every eigenvalue $\lambda \in \C$ with ${\rm Re}(\lambda) \neq 0$,
there is an eigenvalue $-\lambda \in \C$. For instance, if $E''(\phi)$ is invariant under the parity transformation,
then the Hamiltonian symmetry holds and if $v(x)$ is the eigenvector of the spectral
problem (\ref{spectrum}) for $\lambda$, then $v(-x)$ is the eigenvector of the spectral
problem (\ref{spectrum}) for $-\lambda$.

The nonlinear wave $\phi$
is {\em spectrally stable} if $\sigma(\partial_x \mathcal{L}) \subset i \mathbb{R}$
and it is {\em spectrally unstable} if there is $\lambda_0 \in \sigma(\partial_x \mathcal{L})$
such that ${\rm Re}(\lambda_0) > 0$,
where $\sigma(\partial_x \mathcal{L})$ denotes the spectrum of the non-self-adjoint
eigenvalue problem (\ref{spectrum}). The corresponding eigenvector $v$ for an
eigenvalue $\lambda \in \sigma(\partial_x \mathcal{L})$ belongs to the function space
$D(\partial_x \mathcal{L}) \cap \dot{H}^{-1}(\R) \subset L^2(\mathbb{R})$, where $\dot{H}^{-1}(\mathbb{R})$ is the space of
all distributions with square integrable anti-derivatives. In other words, if $v \in \dot{H}^{-1}(\mathbb{R})$, then
$\partial_x^{-1} v \in L^2(\mathbb{R})$.

We assume that the unbounded self-adjoint operator $\mathcal{L}$ from $D(\mathcal{L}) \subset L^2(\R)$
to $L^2(\R)$ is given by the sum of two operators
$\mathcal{L}_0$ and $K_{\mathcal{L}}$, where $\mathcal{L}_0$ is a strongly elliptic unbounded
operator with constant coefficients and $K_{\mathcal{L}}$ is a relatively compact
perturbation of $\mathcal{L}_0$. Using the Fourier transform $\mathcal{F}$ on $L^2(\R)$, we
define the image of $\mathcal{L}_0$ as follows:
$$
\mathcal{F}(\mathcal{L}_0 u)(k) = \hat{\mathcal{L}}_0(k) \mathcal{F}(u)(k), \quad k \in \mathbb{R}.
$$
Since $\mathcal{L}_0$ is unbounded, a coercivity condition holds to yield
$\hat{\mathcal{L}}_0(k) \to \infty$ as $|k| \to \infty$.
We will further assume the following generic assumptions.

\begin{itemize}
\item[(H1)] There is $c_0 > 0$ such that $\hat{\mathcal{L}}_0(k) \geq c_0$ for all $k \in \mathbb{R}$.
By Weyl's theorem, this implies that the essential spectrum of $\mathcal{L}$ (denoted as $\sigma_e(\mathcal{L})$)
is bounded away from zero by a positive number.

\item[(H2)] The discrete spectrum of $\mathcal{L}$ (denoted as $\sigma_d(\mathcal{L})$)
includes a finite number $n(\mathcal{L})$ of negative eigenvalues with eigenvectors in $D(\mathcal{L})$.

\item[(H3)] ${\rm Ker}(\mathcal{L}) = {\rm span}\{f_0\}$ with $f_0 \in D(\mathcal{L}) \cap D(\partial_x \mathcal{L} \partial_x) \cap \dot{H}^{-1}(\mathbb{R})$, so that $\phi_0 = \partial_x^{-1} f_0 \in L^2(\mathbb{R})$\footnote{The restrictive 
    assumption $f_0 \in D(\partial_x \mathcal{L} \partial_x)$ is needed because $f_0$ defines later the projection operator $P$ 
    in the generalized eigenvalue problem (\ref{generalized_eigenvalue_problem}).}.

\item[(H4)] $\langle \mathcal{L}^{-1} \phi_0,\phi_0 \rangle \neq 0$.
The value of $\langle \mathcal{L}^{-1} \phi_0,\phi_0 \rangle$ is finite
because $\langle f_0, \phi_0 \rangle = \langle \partial_x \phi_0, \phi_0 \rangle = 0$.
\end{itemize}

Under these generic assumptions, we obtain the instability index count,
which is analogous to instability index count for NLS-type evolution equations
\cite{ChP10,CPV05,KKS04,P05} (see also Chapter 4 in \cite{Pel-book}).
To formulate the theorem, let us define the following
numbers for the eigenvalue problem (\ref{spectrum}) with the
account of algebraic multiplicity of eigenvalues:
\begin{itemize}
\item $N_r$ is the number of real positive eigenvalues $\lambda$.

\item $N_c$ is the number of complex eigenvalues $\lambda$ in the first open quadrant of $\mathbb{C}$.

\item $N_i^-$ is the total negative Krein index\footnote{The negative Krein index of
an invariant subspace $E_{\lambda} \subset L^2$ of the spectral stability problem (\ref{spectrum})
associated with an eigenvalue $\lambda \in i \mathbb{R}$ is the
number of non-positive eigenvalues of $\langle \mathcal{L}|_{E_{\lambda}} u, u \rangle$.
These eigenvalues of $\langle \mathcal{L}|_{E_{\lambda}} u, u \rangle$
can not be zero if $\lambda$ is an isolated eigenvalue but may include
zero eigenvalue if $\lambda$ is an embedded eigenvalue.}
associated with the number of imaginary (possibly, embedded)
eigenvalues $\lambda$ with ${\rm Im}(\lambda) > 0$.
\end{itemize}
Our main result is the following theorem.

\begin{Theorem}
\label{main-theorem}
Assume (H1)--(H4). Then,
\begin{equation}
\label{closure-relation1} N_r + 2 N_c + 2 N_i^- = n(\mathcal{L}) - n_0,
\end{equation}
where $n_0 = 1$ if $\langle \mathcal{L}^{-1} \phi_0,\phi_0 \rangle < 0$ and $n_0
= 0$ if $\langle \mathcal{L}^{-1} \phi_0,\phi_0 \rangle > 0$.
\end{Theorem}

Section 2 contains historical notes devoted to Theorem \ref{main-theorem} and other recent relevant results.
The proof of Theorem \ref{main-theorem} is developed in Section 3.
A generalization of Theorem \ref{main-theorem} for a periodic nonlinear wave $\phi$
is given in Section 4. Section 5 discusses further possible developments in the area.

\section{Historical remarks and examples}

The Hamiltonian functional $E(u)$ conserves in time $t$ in the KdV-type evolution
equation (\ref{Ham-system}). For many KdV-type evolution equations, there exists typically
another conserved $C^2$ functional $P(u)$, called the momentum functional.
For example, for the general fifth-order KdV equation \cite{CG97,CG94},
\begin{equation}
\label{general-KdV} u_t = a_1 u_x - a_2 u_{xxx} + a_3 u_{xxxxx} + 3
b_1 u u_x - b_2 \left( u u_{xxx} + 2 u_x u_{xx} \right) + 6 b_3 u^2
u_x,
\end{equation}
where $(a_1,a_2,a_3,b_1,b_2,b_3)$ are real,
the energy functional $E(u)$ is well defined in $H^2(\mathbb{R})$,
\begin{equation}
\label{general-Hamiltonian} E(u) = \frac{1}{2} \int_{\mathbb{R}}
\left( a_1 u^2 + a_2 u_x^2 + a_3 u_{xx}^2 + b_1 u^3 + b_2 u u_x^2 +
b_3 u^4 \right) dx,
\end{equation}
whereas the momentum functional is $P(u) = \| u \|^2$. Without
loss of generality, we assume that the phase speed for linear waves
in the fifth-order KdV equation (\ref{general-KdV}) is non-negative.
Using the Fourier transform $\mathcal{F}$, we express this assumption as follows:
\begin{equation}
\label{wave-speed}
c_{\rm wave}(k) = a_1 + a_2 k^2 + a_3 k^4 \geq 0, \quad k \in \mathbb{R}.
\end{equation}
Assumption (\ref{wave-speed}) is needed for assumption (H1) in Theorem \ref{main-theorem}.

Besides the translational parameter $x_0$ in $\phi(x-x_0)$,
the nonlinear wave $\phi$ has typically another free parameter $c$
for the constant speed. With the account of speed $c$, the nonlinear wave
is a critical point of the {\em extended} energy functional
$E_c(u) := E(u) + c P(u)$. The second variation of $E_c$ defines
the self-adjoint operator $\mathcal{L}_c := E''(\phi) + c P''(\phi)$.

General stability-instability results for the critical points of
$E_c(u)$ were obtained in \cite{BSS87,GSS,SS90}, based on the
assumption that the self-adjoint linearized operator $\mathcal{L}_c$ has
exactly one negative eigenvalue and a simple zero eigenvalue.
By a different method involving modulation equations, Lyapunov stability
of positive travelling waves $\phi$ was also proved by Weinstein \cite{Wen}.
In the consequent two influential papers, Pego and Weinstein developed
Evans function analysis of spectral stability \cite{PW92} and analysis
of asymptotic stability in exponentially weighted spaces \cite{PW94} in the context of
a generalized KdV equation.

These general results correspond to the case $n(\mathcal{L}) = 1$ in Theorem \ref{main-theorem}
(see also \cite{Pava} for a recollection of these and many other results). More recently, questions have
been raised on spectral stability of KdV-type nonlinear waves in the cases
where $n(\mathcal{L}) > 1$, which are known for equations of the integrable
KdV hierarchy \cite{MS,KP}. The result of Theorem \ref{main-theorem} was
already claimed as early as 2006 in the context of the fifth-order KdV
equation (\ref{general-KdV}) \cite{ChPel-preprint}, although the final
version of this paper was published without the example of
the fifth-order KdV equation \cite{ChP10}.
%\footnote{Peer-reviewing
%policies in some American journals resemble now the worst traditions of censorship in former Soviet Union.
%Referee reports are biased towards ``prevailing theories", editors are sceptical
%against submissions from outside of the ``inner circle" of their colleagues or friends,
%and blind rejections of good papers with new results are not uncommon.}.
Since that time, a weaker result was obtained by Lin \cite{Lin}
and a nearly identical result was outlined very recently
by Kapitula \& Stefanov \cite{KapStef}. Periodic waves of KdV-type nonlinear
evolution equations were treated in \cite{Natali,Angulo,BronKap,DecKap,KapHar}, where results
similar to Theorem \ref{main-theorem} were obtained. Therefore, it makes fair to
restore the original proof of Theorem \ref{main-theorem} following the lines
of \cite{ChPel-preprint} and to show how naturally the instability index count for
both solitary and periodic waves can be
adopted from a general theory in Pontryagin's space \cite{Pontryagin}.
This task is achieved in the present paper with the main goal to show the universality
and simplicity of the proof of the instability index counts by using
the generalized eigenvalue problem (that is, the linear operator pencil
in the terminology of the recent review \cite{KollarMiller}).

In the context of the general fifth-order KdV equation (\ref{general-KdV}),
specific studies of Lyapunov stability of travelling solitary waves
were reported in \cite{DK99,IS92} with the energy-momentum methods.
In particular, since the solitary wave satisfies the fourth-order differential equation,
\begin{equation}
\label{ODE-fifth-order} a_3 \phi'''' - a_2 \phi'' + (a_1 + c) \phi +
\frac{3}{2} b_1 \phi^2 - \frac{1}{2} b_2 \left( 2 \phi \phi'' +
(\phi')^2 \right) + 2 b_3 \phi^3 = 0,
\end{equation}
one can verify by direct computations that $\mathcal{L}_c \phi' = 0$ and $\mathcal{L}_c \partial_c \phi = -\phi$, where
the prime denotes differentiation in $x$ and $\partial_c$ denotes differentiation in $c$,
whereas
\begin{equation}
\label{operator-L} \mathcal{L}_c := a_3 \frac{d^4}{d x^4} - a_2
\frac{d^2}{d x^2} + a_1 + c + 3 b_1 \phi(x) - b_2 \frac{d}{d x}
\phi(x) \frac{d}{d x} - b_2 \phi''(x) + 6 b_3 \phi^2(x).
\end{equation}
Assuming existence and uniqueness (up to translational invariance)
of a solitary wave $\phi \in H^2(\R)$ with the exponential decay
at infinity for $c > 0$ (see \cite{CG97,K97,L99} for existence results),
we realize that the operator $\mathcal{L}_c$ satisfies assumptions (H1)--(H4) of Theorem
\ref{main-theorem} with $\hat{\mathcal{L}}_0(k) = c + c_{\rm wave}(k) \geq c > 0$,
$\phi_0 = \phi$, and
$$
\langle \mathcal{L}_c^{-1} \phi,\phi \rangle = - \langle \partial_c \phi, \phi \rangle =
-\frac{1}{2} \frac{d}{dc} \| \phi \|^2.
$$
If $n(\mathcal{L}_c) = 1$, the result of Theorem \ref{main-theorem} gives stability
of a solitary wave if $\frac{d}{dc} \| \phi \|^2 > 0$ and instability
if $\frac{d}{dc} \| \phi \|^2 < 0$, which coincides with the results of
the orbital stability theory \cite{Pava,BSS87,Wen}.

Spectral stability of one-humped solitary waves
in the fifth-order KdV equation was studied numerically in \cite{BDG03},
with the use of the symplectic Evans matrix \cite{BD02}. Because
$n(\mathcal{L}_c) = 1$ and $\frac{d}{dc} \| \phi \|^2 > 0$ were found,
the one-humped solitary waves were shown to be spectrally stable.

One-humped and two-humped solitary waves in the fifth-order KdV equation
were numerically approximated in \cite{ChP06} with a spectral method.
Numerical results on eigenvalues of the spectral
problem (\ref{spectrum}) were found in full correspondence with the result of Theorem \ref{main-theorem}.
In particular, two-humped solutions have either $n(\mathcal{L}_c) = 2$ or $n(\mathcal{L}_c) = 3$,
depending whether the individual solitary waves
form a bound state at the non-degenerate minimum or maximum points
of the effective interaction potential. Since $\frac{d}{dc} \| \phi \|^2 > 0$
for all these solitary waves, the two-humped solutions with $n(\mathcal{L}_c) = 2$ are
unstable with $N_r = 1$. Nevertheless, the two-humped solutions with $n(\mathcal{L}_c) = 3$
are spectrally stable, because the single pair of embedded eigenvalues with
negative Krein signature $N_i^- = 1$ is structurally stable with respect to
parameter continuations \cite{ChP06}. Similar results were also observed
numerically with the computations of the Maslov index for solitary waves
in the fifth-order KdV equation \cite{Chardard}.

To finish these remarks, we also mention a similar instability index count
obtained for dark solitons in the defocusing NLS equation with an external potential
\cite{PelKev}. Although the symplectic operator for the NLS equation is invertible,
the spectral stability problem for dark solitons (solitary waves with nonzero boundary
conditions) is defined in terms of a linear self-adjoint operator, where the positive
essential spectrum touches zero. Nevertheless, the theory
from \cite{ChP10} was successfully applied to the count of unstable eigenvalues for a dark
soliton in a spatially localized potential and illustrated with a number of prototypical
examples in \cite{PelKev}. In this context, a dark soliton
persists in a small localized potential if it is located at the non-degenerate minimum or maximum points of the effective
potential and is spectrally unstable in both cases. At the maximum point, the dark soliton
is unstable with one real eigenvalue $N_r = 1$, whereas at the minimum point, it is
unstable with two complex eigenvalues $N_c = 1$. The embedded imaginary eigenvalues with negative Krein signature
are structurally unstable with respect to parameter continuations in the defocusing NLS equations
and bifurcate into complex unstable eigenvalues (see \cite{PelKev} for precise computations
of these instabilities by using the Evans function for dark solitons).

\section{Proof of Theorem \ref{main-theorem}}

We consider the spectral problem (\ref{spectrum}), where the self-adjoint
operator $\mathcal{L} = E''(\phi)$ satisfies the assumptions (H1)--(H4) of Theorem \ref{main-theorem}.

The proof of the standard instability index count
\cite{ChP10,CPV05,KKS04,P05} is not applicable to the KdV-type
evolution equations because the symplectic operator $\partial_x$ is
not invertible. Nevertheless, the range of the self-adjoint operator $\mathcal{L}$
is defined in $L^2(\mathbb{R})$, hence bootstrapping arguments imply that the
eigenvector $v$ of the spectral problem (\ref{spectrum}) with $\lambda \neq 0$
belongs to the function space $D(\partial_x \mathcal{L}) \cap \dot{H}^{-1}(\mathbb{R}) \subset L^2(\R)$.
Therefore, we can define $w = \partial_x^{-1} v \in L^2(\mathbb{R})$ and formally
extend the spectral problem (\ref{spectrum}) to the system of two coupled equations:
\begin{equation}
\label{coupled-problem} \mathcal{M} w = -\lambda v, \quad \mathcal{L} v = \lambda w,
\end{equation}
where $\mathcal{M} := - \partial_x \mathcal{L} \partial_x$,
$v \in D(\partial_x \mathcal{L}) \cap \dot{H}^{-1}(\R) \subset L^2(\R)$,
and $w \in D(\partial_x \mathcal{L} \partial_x) \subset L^2(\R)$. The coupled system
(\ref{coupled-problem}) is equivalent to the squared eigenvalue
problem $\partial_x \mathcal{L} \partial_x \mathcal{L} v = \lambda^2 v$.
We show now that if the coupled system (\ref{coupled-problem}) has an eigenvalue $\lambda_0 \neq 0$,
then it has another eigenvalue $-\lambda_0$ and these two eigenvalues are equivalent to
the pair of eigenvalues $\lambda_0$ and $-\lambda_0$ of the spectral problem
(\ref{spectrum}). For simplicity of presentation, we only consider the case of simple nonzero
eigenvalues in this work.

\begin{Proposition}
The coupled system (\ref{coupled-problem}) has a pair of simple eigenvalues
$\pm \lambda_0 \neq 0$ with the eigenvectors $(v_0,\pm w_0) \in D(\partial_x \mathcal{L}) \cap \dot{H}^{-1}(\R) \times
D(\partial_x \mathcal{L} \partial_x)$ if and only if the spectral problem (\ref{spectrum}) has a pair of simple
eigenvalues $\pm \lambda_0$ with the eigenvectors $v_{\pm} = v_0 \pm \partial_x w_0 \in D(\partial_x \mathcal{L}) \cap \dot{H}^{-1}(\R)$.
\label{equivalence-lemma}
\end{Proposition}

\begin{Proof}
By the symmetry, if $\lambda_0 \neq 0$ is a simple eigenvalue of
the coupled system (\ref{coupled-problem}) with
the eigenvector $(v_0,w_0) \in D(\partial_x \mathcal{L}) \cap \dot{H}^{-1}(\R) \times
D(\partial_x \mathcal{L} \partial_x)$, then $-\lambda_0$ is also a simple eigenvalue
of the coupled system (\ref{coupled-problem}) with the eigenvector $(v_0,-w_0)$.
Moreover, $v_0$ and $w_0$ are linearly independent.

We differentiate the second equation of the coupled system (\ref{coupled-problem}) for
the eigenvalue $\lambda_0$ and add or subtract the first equation of the system to obtain
$$
\partial_x \mathcal{L} (v_0 \pm \partial_x w_0) = \pm \lambda_0 (v_0 \pm \partial_x w_0).
$$
Therefore, $v_0 \pm \partial_x w_0 \in {\rm Ker}(\partial_x \mathcal{L} \mp \lambda_0)$.
By the Hamiltonian symmetry, if $\lambda_0 \in \sigma(\partial_x \mathcal{L})$, then
$-\lambda_0 \in \sigma(\partial_x \mathcal{L})$, whereas the algebraic multiplicity
of eigenvalues in $\sigma(\partial_x \mathcal{L} \partial_x \mathcal{L})$
equals the algebraic multiplicity of eigenvalues in the coupled system (\ref{coupled-problem}).
This guarantees that the two eigenvectors $(v_0,\pm w_0)$ of the coupled system (\ref{coupled-problem})
generate two linearly independent eigenvectors  $v_{\pm} = v_0 \pm \partial_x w_0$ for $\pm \lambda_0
\in \sigma(\partial_x \mathcal{L})$.

To check the converse statement, we assume that $v_{\pm}$ are linearly independent eigenvectors of
$\partial_x \mathcal{L}$ for the eigenvalues $\pm \lambda_0$ in $D(\partial_x \mathcal{L}) \cap \dot{H}^{-1}(\R)$.
Then, we define nonzero functions
\begin{equation}
\label{nonzero-functions}
v_0 := \frac{1}{2} (v_+ + v_-), \quad w_0 := \frac{1}{2}( \partial_x^{-1} v_+ - \partial_x^{-1} v_-),
\end{equation}
and obtain
$$
\partial_x \mathcal{L} v_0 = \frac{1}{2} \partial_x \mathcal{L} (v_+ + v_-) = \frac{1}{2} \lambda_0 ( v_+ - v_-) = \lambda_0 \partial_x w_0,
$$
so that the integration gives $\mathcal{L} v_0 = \lambda_0 w_0$, that is, the second equation of
the coupled system (\ref{coupled-problem}).
Similarly, we check the first equation of the coupled system (\ref{coupled-problem}) with $\mathcal{M} w_0 = -\lambda_0 v_0$.
Therefore, the two eigenvectors $v_{\pm}$ for $\pm \lambda_0 \in \sigma(\partial_x \mathcal{L})$
generate two linearly independent eigenvectors $(v_0,\pm w_0)$ of the coupled system (\ref{coupled-problem})
for a pair of eigenvalues $\pm \lambda_0 \neq 0$.
\end{Proof}

\begin{Remark}
In many KdV-type evolution equations including the fifth-order KdV equation (\ref{general-KdV}),
the Hamiltonian symmetry of the spectral problem (\ref{spectrum}) follows from the parity transformation
of the eigenvectors, if the nonlinear wave $\phi$ is symmetric with respect to $x$. Therefore, if
$v_+(x)$ is a solution of $\partial_x \mathcal{L} v_+ = \lambda_0 v_+$, then $v_-(x) := v_+(-x)$ is
a solution of $\partial_x \mathcal{L} v_- = -\lambda_0 v_-$. Under this transformation,
the components $v_0$ and $w_0$ of the coupled system (\ref{coupled-problem})
for a simple eigenvalue $\lambda_0$ are either even or odd functions
with respect to $x$, whereas $v_{\pm}$ are neither even nor odd.
\end{Remark}

To study the spectrum of the coupled system (\ref{coupled-problem}), we shall first understand
the spectrum of operator $\mathcal{M}$. Recall again that $\sigma_e$ and $\sigma_d$ denote
the essential and discrete spectra. From assumptions
(H1)--(H4), we obtain the following properties of operator $\mathcal{M}$.

\begin{Lemma}
\label{lemma-operator-M}
Under assumptions (H1)--(H4) on $\mathcal{L}$, operator $\mathcal{M}$ can be extended to a self-adjoint
operator with a dense domain $D(\mathcal{M})$ in $L^2(\mathbb{R})$ satisfying the following properties:
\begin{itemize}
\item[{\rm (H1$^{\prime}$)}] $\sigma_e(\mathcal{M}) \geq 0$.

\item[{\rm (H2$^{\prime}$)}] $\sigma_d(\mathcal{M})$ includes $n(\mathcal{L})$ negative eigenvalues with eigenvectors in $D(\mathcal{M})$.

\item[{\rm (H3$^{\prime}$)}] ${\rm Ker}(\mathcal{M}) = {\rm span}\{\phi_0\}$.

\item[{\rm (H4$^{\prime}$)}] $\langle \mathcal{M}^{-1} f_0,f_0 \rangle$ is finite and nonzero.
\end{itemize}
\end{Lemma}

\begin{Proof}
From the decomposition $\mathcal{L} = \mathcal{L}_0 + K_{\mathcal{L}}$,
we have the decomposition $\mathcal{M} = \mathcal{M}_0 + K_{\mathcal{M}}$, where $K_{\mathcal{M}} = -\partial_x K_{\mathcal{L}} \partial_x$
is a relatively compact perturbation of $\mathcal{M}_0 = -\partial_x \mathcal{L}_0 \partial_x$. Since $\mathcal{M}_0$ is
a linear operator with constant coefficients, we use the Fourier transform $\mathcal{F}$ on
$L^2(\R)$ to find the image of $\mathcal{M}_0$ as follows:
$$
\hat{\mathcal{M}}_0(k) = k^2 \hat{\mathcal{L}}_0(k) \geq 0 \quad \mbox{\rm for all } \; k \in \mathbb{R}.
$$
Since $\hat{\mathcal{L}}_0(k) \geq c_0$ by assumption (H1), we have $\hat{\mathcal{M}}_0(k) \geq 0$ for all $k \in \mathbb{R}$.
By Weyl's theorem, this implies that $\sigma_e(\mathcal{M})$ is non-negative, that is, (H1$^{\prime}$) holds.

Since ${\rm Ker}(\partial_x) = {\rm span}\{0\}$ in $L^2(\R)$, (H3$^{\prime}$)
follows from (H3) by direct computations:
$$
-\partial_x \mathcal{L} \partial_x f = 0 \;\; \Rightarrow \;\; \mathcal{L} \partial_x f = 0 \;\; \Rightarrow \;\;
\partial_x f \in {\rm span}\{\partial_x \phi_0\}, \;\; \Rightarrow \;\; f \in {\rm span}\{\phi_0\}.
$$
Furthermore, $\langle \mathcal{M}^{-1} f_0, f_0 \rangle$ exists because $\partial_x^{-1}$ is well defined
in $L^2(\R)$ on functions in $L^2(\R) \cap \dot{H}^{-1}(\R)$. As a result,
(H4$^{\prime}$) follows from (H4) by means of integration by parts:
$$
\langle \mathcal{M}^{-1} f_0, f_0 \rangle = \langle \mathcal{M}^{-1} \partial_x \phi_0, \partial_x \phi_0 \rangle =
- \langle \partial_x \mathcal{M}^{-1} \partial_x \phi_0,\phi_0 \rangle = \langle \mathcal{L}^{-1} \phi_0, \phi_0 \rangle.
$$

It remains to prove (H2$^{\prime}$). The negative eigenvalues of $\sigma_d(\mathcal{M})$
are defined from the eigenvalue problem $\mathcal{M} f = \lambda f$, which is rewritten
in the following form:
\begin{equation}
\label{H2-1}
-\partial_x \mathcal{L} \partial_x f = \lambda f,
\quad f \in D(\partial_x \mathcal{L} \partial_x) \cap \dot{H}^{-1}(\R) \subset L^2(\R).
\end{equation}
Since $f \in \dot{H}^{-1}(\R)$, there exists $g = \partial_x^{-1} f \in L^2(\mathbb{R})$
such that the spectral problem (\ref{H2-1}) can be written in the equivalent form
\begin{equation}
\label{H2-2}
-\mathcal{L} \partial_x^2 g = \lambda g, \quad g \in D(\mathcal{L} \partial_x^2) \subset L^2(\R).
\end{equation}
For any $\epsilon > 0$, the positive operator $(\epsilon - \partial_x^2)$ is invertible
and the inverse operator is defined by the integral representation
\begin{equation}
\label{integral-representation}
(\epsilon - \partial_x^2)^{-1} f(x) = \frac{1}{2 \epsilon^{1/2}}
\int_{\mathbb{R}} e^{-\epsilon^{1/2} |x-y|} f(y) dy, \quad f \in L^2(\mathbb{R}).
\end{equation}
Using this integral representation, we define a smoothen version
of the eigenvalue problem (\ref{H2-2}) for $h = (\epsilon - \partial_x^2) g$:
\begin{equation}
\label{H2-3}
\mathcal{L} h = \lambda (\epsilon - \partial_x^2)^{-1} h, \quad h \in D(\mathcal{L}) \subset L^2(\R).
\end{equation}
Since $(\epsilon - \partial_x^2)^{-1}$ is a positive bounded self-adjoint operator for any $\epsilon > 0$,
Sylvester's Law of Inertia (Theorem 4.2 in \cite{Pel-book}) applies and the number of
negative and zero eigenvalues of the spectral problem (\ref{H2-3}) corresponds to the number
of negative and zero eigenvalues of the operator $\mathcal{L}$. By assumptions
(H2) and (H3), there are exactly $n(\mathcal{L})$ negative eigenvalues of
the eigenvalue problem (\ref{H2-3}) and a simple zero eigenvalue for any $\epsilon > 0$.

Because the integral representation (\ref{integral-representation}) diverges as $\epsilon \downarrow 0$
and the operator $\mathcal{L}$ is bounded from below,
the negative eigenvalues of  the spectral problem (\ref{H2-3}) for $\epsilon > 0$ are bounded from
below but may a priori approach to zero as $\epsilon \downarrow 0$. However,
since the kernel of $(\epsilon -  \partial_x^2) \mathcal{L}$ is simple for
any $\epsilon \geq 0$, the negative eigenvalues are bounded away from zero as $\epsilon \downarrow 0$.
%On the other hand, by bootstrapping arguments,
%if $h \in D(\mathcal{L})$, then $g = (\epsilon - \partial_x^2)^{-1} h \in
%D(\mathcal{L} \partial_x^2) \subset L^2(\R)$
%is an eigenvector of $\mathcal{L}(\epsilon - \partial_x^2) g = \lambda g$
%for any $\epsilon \geq 0$. The negative eigenvalues of the latter problem
%are semi-continuous (remains finite and negative) as $\epsilon \downarrow 0$
As a result, the spectral problem (\ref{H2-2}) also has $n(\mathcal{L})$ negative eigenvalues, that is, (H2$^{\prime}$) holds\footnote{Another smoothen version of the same eigenvalue problem (\ref{H2-2}) is
$(\epsilon - \partial_x^2) g = \lambda \mathcal{L}^{-1} g$ for all $g \in H^2(\R) \cap [{\rm span}\{f_0\}]^{\perp}$.
By the same Sylvester's Law of Inertia, there are exactly $n(\mathcal{L})$ negative eigenvalues of
this eigenvalue problem for any $\epsilon > 0$ and the bootstrapping arguments give
$g \in D(\mathcal{L} \partial_x^2) \subset L^2(\R)$ for the corresponding eigenvectors.}.
\end{Proof}

We convert now the coupled system (\ref{coupled-problem}) for a pair of simple eigenvalues $\pm \lambda \neq 0$ to a generalized eigenvalue
problem for a double eigenvalue. By assumptions (H1) and (H3), the zero eigenvalue of $\mathcal{L}$
is bounded away from the essential spectrum of $\mathcal{L}$.
Let $P$ be the orthogonal projection from $L^2(\R)$ to $[{\rm span}\{f_0\}]^{\perp} \subset L^2(\R)$.
The following result establishes this equivalence.

\begin{Proposition}
\label{proposition-generalized-eigenvalue}
The coupled system (\ref{coupled-problem}) has a pair of simple eigenvalues
$\pm \lambda \neq 0$ with eigenvectors $(v,\pm w)$ if and only if the generalized
eigenvalue problem
\begin{equation}
\label{generalized_eigenvalue_problem} A w = \gamma K w, \quad w \in
\mathcal{H} := D(\mathcal{M}) \cap  [{\rm span}\{f_0\}]^{\perp} \subset L^2(\R),
\end{equation}
where $A := P \mathcal{M} P$ and $K := P \mathcal{L}^{-1} P$,
has a double eigenvalue $\gamma = - \lambda^2 \neq 0$ with linearly
independent eigenvectors $\partial_x^{-1} v$ and $w$.
\end{Proposition}

\begin{Proof}
Let $\pm \lambda \neq 0$ be a pair of simple eigenvalues of the coupled system  (\ref{coupled-problem})
with the eigenvectors $(v,\pm w) \in D(\partial_x \mathcal{L}) \cap \dot{H}^{-1}(\R) \times D(\partial_x \mathcal{L} \partial_x)$.
Because $\lambda \neq 0$, we have $w = Pw$, that is, $w$ is in the range of $\mathcal{L}$.
As a result, the second equation of the coupled
system (\ref{coupled-problem}) can be written in the equivalent form:
\begin{equation}
\label{equation-w}
v = \lambda P \mathcal{L}^{-1} P w + v_0, \quad v_0 \in {\rm Ker}(\mathcal{L}).
\end{equation}
Substituting $v$ into the first equation of the coupled system (\ref{coupled-problem})
and using the projection operator $P$ again, we obtain a closed equation for $w$
\begin{equation}
\label{equation-u}
P \mathcal{M} P w = - \lambda^2 P \mathcal{L}^{-1} P w, \quad w \in D(\mathcal{M}) \cap [{\rm span}\{f_0\}]^{\perp} \subset L^2(\R)
\end{equation}
and a unique expression for $v_0$
\begin{equation}
\label{equation-v}
v_0 = -\frac{1}{\lambda} (I - P) \mathcal{M} P w,
\end{equation}
where $\lambda \neq 0$ and $(I-P)$ is the orthogonal projection from
$L^2(\R)$ to ${\rm Ker}(\mathcal{L})$. Therefore, it follows
from equation (\ref{equation-u}) that $\gamma = -\lambda^2$ is
an eigenvalue of the generalized eigenvalue problem
(\ref{generalized_eigenvalue_problem}) with an eigenvector $w$.

To show that this $\gamma$ is a double
eigenvalue of the generalized eigenvalue problem
(\ref{generalized_eigenvalue_problem}), we note that the coupled system
(\ref{coupled-problem}) is invariant with respect to
the transformation
\begin{equation}
\label{transformation-formula}
\partial_x w \to v \quad \mbox{\rm and} \quad \partial_x^{-1} v \to w.
\end{equation}
Therefore, $\partial_x^{-1} v$ is another eigenvector of the generalized eigenvalue problem
(\ref{generalized_eigenvalue_problem}) for the same $\gamma$.
By Proposition \ref{equivalence-lemma}, see the equivalence formula (\ref{nonzero-functions}),
$\partial_x^{-1} v$ and $w$ are linearly independent.

In the opposite direction, if $\gamma \neq 0$ is a double eigenvalue of the generalized
eigenvalue problem (\ref{generalized_eigenvalue_problem})
with linearly independent eigenvectors $w_1$ and $w_2$,
then for each eigenvector $w_1$ or $w_2$, we define $v_0$ by (\ref{equation-v})
and $v$ by (\ref{equation-w}), which yields linearly independent
components $v_1$ and $v_2$ of the coupled system (\ref{coupled-problem})
for the eigenvalue $\lambda = (-\gamma)^{1/2}$. The eigenvalue $\lambda = (-\gamma)^{1/2}$
must be simple (or the multiplicity of the eigenvalue $\gamma$ in the generalized eigenvalue problem
(\ref{generalized_eigenvalue_problem}) exceeds two), therefore, the
transformation (\ref{transformation-formula}) yields the correspondence
$v_1 = \partial_x w_2$ and $v_2 = \partial_x w_1$. In other words, only $(v_1,w_1)$
is a linearly independent eigenvector of the coupled system (\ref{coupled-problem}) for
the simple eigenvalue $\lambda = (-\gamma)^{1/2}$. The other simple eigenvalue
$-\lambda = -(-\gamma)^{1/2}$ exists by the symmetry of the coupled system (\ref{coupled-problem})
with the eigenvector $(v_1,-w_1)$.
\end{Proof}

The generalized eigenvalue
problem (\ref{generalized_eigenvalue_problem})
for unbounded self-adjoint differential operators $A$ and $K$
with strictly positive essential spectrum was studied by
Chugunova \& Pelinovsky \cite{ChP10} in
Pontryagin's space \cite{Pontryagin}. Here we report
the modification of the analysis needed to treat the case when
the bottom of the essential spectrum of $A$ touches zero.
We shall first prove that a deformation of $A$ to $A_{\delta} := A + \delta K$ for
a small positive number $\delta$ shifts the essential spectrum away from zero.

\begin{Lemma}
For small positive values of $\delta$, there is a positive $\delta$-independent constant $d_0$ such that
\begin{equation}
\sigma_e(A_{\delta}) \geq d_0 \delta.
\end{equation}
\label{lemma-ess-spectrum}
\end{Lemma}

\begin{Proof}
Since $\mathcal{M}$ and $\mathcal{L}$ are represented by the
relatively compact perturbations of operators $\mathcal{M}_0$ and $\mathcal{L}_0$ with constant coefficients,
we can use the Fourier transform $\mathcal{F}$ on $L^2(\R)$ to compute:
$$
\mathcal{F}(\mathcal{M}_0 + \delta \mathcal{L}_0^{-1}) = k^2 \hat{\mathcal{L}}_0(k)
+ \delta \hat{\mathcal{L}}_0^{-1}(k), \quad k \in \mathbb{R},
$$
where $\hat{\mathcal{L}}_0(k) \geq c_0 > 0$ for some $c_0$ by (H1).

Let $k_{\delta}$ denote the positive
global minimum of this function. By coercivity of $k^2 \hat{\mathcal{L}}_0(k)$, the global minimum
is achieved at a finite value of $k$ for small positive values of $\delta$
and there is a $\delta$-independent
positive constant $K_0$ such that $k_{\delta} \in [0,K_0]$. But then, there is
a $\delta$-independent positive constant $d_0$ such that
$\hat{\mathcal{L}}_0^{-1}(k_{\delta}) \geq d_0$ and
$k^2 \hat{\mathcal{L}}_0(k) + \delta \hat{\mathcal{L}}_0^{-1}(k) \geq d_0 \delta$
for small positive $\delta$.
\end{Proof}

By Lemma \ref{lemma-ess-spectrum}, the essential spectrum of $A_{\delta}$ for a small positive $\delta$
is strictly positive. Also, the kernel of $A_{\delta}$ is empty for a small positive $\delta$
because if $f \in {\rm Ker}(A_{\delta})$, then $A f = -\delta K f$ but the negative eigenvalues do not accumulate
near zero, thanks to the decomposition $\mathcal{L} = \mathcal{L}_0 + K_{\mathcal{L}}$ with
a relatively compact perturbation $K_{\mathcal{L}}$. Therefore,
there exists a small positive number $\delta$ such
that operator $A_{\delta}$ is continuously invertible in ${\cal H}$ and
the generalized eigenvalue problem
(\ref{generalized_eigenvalue_problem}) is rewritten in the shifted
form,
\begin{equation}
\label{shifted_eigenvalue_problem} (A + \delta K) w = (\gamma +
\delta)K w, \qquad w \in {\cal H}.
\end{equation}

By the spectral theory of self-adjoint operators, the Hilbert space ${\cal
H}$ can be equivalently decomposed into two orthogonal sums of
subspaces which are invariant with respect to the operators $K$ and
$A_{\delta}$ for small positive values of $\delta$:
\begin{eqnarray}
\label{decomposition-K} \mathcal{H} = \mathcal{H}_{K}^- \oplus
\mathcal{H}_{K}^+ = \mathcal{H}_{A_{\delta}}^- \oplus \mathcal{H}_{A_{\delta}}^+,
\end{eqnarray}
where notation $-(+)$ stands for invariant subspaces of these operators
related to the negative (positive) spectrum.

Since $P$ is a projection defined by the eigenspace of $\mathcal{L}$ and $K = P \mathcal{L}^{-1} P$, it is obvious that
\begin{equation}
\label{count-K}
{\rm dim}(\mathcal{H}_{K}^-) = n(\mathcal{L}).
\end{equation}
On the other hand,
the number of negative eigenvalues of $A = P\mathcal{M}P$ is related to the
number of negative eigenvalues of $\mathcal{M}$. Compared to the standard
count of negative eigenvalues in constrained Hilbert spaces (Theorem 4.1 in \cite{Pel-book}), the
complication here is that the zero eigenvalue of $\mathcal{M}$ is embedded to the edge of
the essential spectrum of $\mathcal{M}$. In addition, the zero eigenvalue of $A$ is shifted under the
perturbation $\delta K$ in the operator $A_{\delta} = A + \delta K$.
The following two lemmas give the count of negative eigenvalues of $A$ denoted as $n(A)$ and
the count of ${\rm dim}(\mathcal{H}_{A_{\delta}}^-)$ for a small positive number $\delta$.

\begin{Lemma}
\label{lemma-negative-A}
Under assumptions (H1)--(H4) on $\mathcal{L}$, we have
\begin{equation}
n(A) = n(P\mathcal{M}P) = n(\mathcal{M})- n_0 = n(\mathcal{L}) - n_0,
\end{equation}
where $n_0 = 1$ if $\langle \mathcal{L}^{-1} \phi_0,\phi_0 \rangle < 0$ and $n_0
= 0$ if $\langle \mathcal{L}^{-1} \phi_0,\phi_0 \rangle > 0$\footnote{By the standard
technique (Theorem 4.1 in \cite{Pel-book}), we also have $n(\tilde{P}\mathcal{L}\tilde{P}) = n(\mathcal{L}) - n_0$, where $\tilde{P}$ is
an orthogonal projection from $L^2(\R)$ to $[{\rm span}\{\phi_0\}]^{\perp} \subset L^2(\R)$.}.
\end{Lemma}

\begin{Proof}
We study the behavior of the function $F(\mu) = \langle (\mu - \mathcal{M})^{-1} f_0,f_0 \rangle$,
which is well-defined for all $\mu \in \R_- \backslash \sigma(\mathcal{M})$. By (H4$^{\prime}$), it has
the limit as $\mu$ approaches zero from below:
$$
\lim_{\mu \uparrow 0} F(\mu) = - \langle \mathcal{M}^{-1} f_0,f_0 \rangle
= -\langle \mathcal{L}^{-1} \phi_0,\phi_0\rangle \neq 0.
$$
Hence, the assertion of the lemma holds by the standard proof of
Theorem 4.1 in \cite{Pel-book} (where it is formulated and proved in a more general setting).
\end{Proof}

\begin{Lemma}
\label{lemma-A-delta-K}
Under assumptions (H1)--(H4) on $\mathcal{L}$, for a small positive number $\delta$,
we have
\begin{equation}
\label{count-A-delta}
{\rm dim}(\mathcal{H}_{A_{\delta}}^-) = n(\mathcal{L}).
\end{equation}
\end{Lemma}

\begin{Proof}
Negative eigenvalues of $\sigma_d(A_{\delta})$ are defined from the eigenvalue problem:
\begin{equation}
A f + \delta K f = \lambda f \quad f \in \mathcal{H}.
\end{equation}
We use the assumptions that ${\rm Ker}(A) = {\rm span}\{\phi_0\}$ and
$\langle K \phi_0,\phi_0 \rangle = \langle \mathcal{L}^{-1} \phi_0,\phi_0 \rangle \neq 0$\footnote{Note that $P \phi_0 = \phi_0$
because $\langle f_0, \phi_0 \rangle = 0$ and $P$ is an orthogonal projection to $[{\rm span}\{f_0\}]^{\perp}$.}.
Since negative eigenvalues of $A$ in Lemma \ref{lemma-negative-A}
are bounded away from zero, they persist for small positive values of $\delta$.
Let $\mathcal{H}_{\delta}$ denote the orthogonal complement of
the subspace spanned by $n(A)$ eigenvectors corresponding to these
negative eigenvalues of $A + \delta K$ for small positive values of $\delta$.

At $\delta = 0$, we have $\phi_0 \in \mathcal{H}_{\delta = 0}$.
If $\langle K \phi_0,\phi_0 \rangle < 0$, then $A_{\delta} = A + \delta K$ is not positive
definite on $\mathcal{H}_{\delta}$ for small positive $\delta$. Therefore, there is
at least one negative (isolated) eigenvalue of $A_{\delta}$, which
becomes the zero eigenvalue of $A$ as $\delta \to 0$
(the zero eigenvalue of $A$ is embedded at the edge of $\sigma_e(A)$).
Moreover, this is the only small negative eigenvalue of $A_{\delta}$ for small
positive $\delta$\footnote{The edge of $\sigma_e(A)$ may generate additional eigenvalues
by means of edge bifurcations \cite{KapSan}. All these eigenvalues are strictly
positive because they detach from the bottom of $\sigma_e(A_{\delta})$ which is as small
as $\mathcal{O}(\delta)$, whereas the distance of these eigenvalues from the bottom
of $\sigma_e(A_{\delta})$ may only change as a superlinear function of $\delta$ as $\delta \to 0$ \cite{KapSan}.
Therefore, all these eigenvalues via edge bifurcations are necessarily positive.}.
Thus, we conclude that if $\langle K \phi_0,\phi_0 \rangle = \langle \mathcal{L}^{-1} \phi_0, \phi_0 \rangle < 0$, then
$$
{\rm dim}(\mathcal{H}_{A_{\delta}}^-) = n(A) + 1 = n(\mathcal{L}).
$$

On the other hand, if $\langle K \phi_0,\phi_0 \rangle > 0$, the operator $A_{\delta} = A + \delta K$ is
strictly positive on the subspace $\mathcal{H}_{\delta}$ for small positive $\delta$. Therefore,
in this case, we have
$$
{\rm dim}(\mathcal{H}_{A_{\delta}}^-) = n(A) = n(\mathcal{L}).
$$
The assertion of the lemma is proven in both the cases.
\end{Proof}

We are now ready to use Theorem 1 from \cite{ChP10}. Note that although
the theorem was proven under the assumption that the essential spectrum of
$A$ is bounded away from zero, the shift of $A$ to $A_{\delta}$
satisfying $\sigma_e(A_{\delta}) \geq d_0 \delta > 0$ justifies the technique behind
the proof of Theorem 1 in \cite{ChP10} for a small positive number $\delta$.
To formulate the theorem, we introduce some notations
for the numbers of particular eigenvalues $\gamma$
of the generalized eigenvalue problem (\ref{generalized_eigenvalue_problem})
with the account of their algebraic multiplicities.

\begin{itemize}
\item $N_p^-$ ($N_n^-$) is the number of negative eigenvalues $\gamma$
whose (generalized) eigenvectors are associated to the non-negative (non-positive) values of the
quadratic form $\langle K \cdot,\cdot \rangle$.

\item $N_p^+$ ($N_n^+$) is the number of positive eigenvalues $\gamma$
whose (generalized) eigenvectors are associated to the non-negative (non-positive) values of the
quadratic form $\langle K \cdot,\cdot \rangle$.

\item $N_p^0$ ($N_n^0$) is the multiplicity of zero eigenvalue whose
(generalized) eigenvectors are associated to the non-negative (non-positive) values of the
quadratic form $\langle K \cdot,\cdot \rangle$.

\item $N_{c^+}$ ($N_{c^-}$) is the number of complex eigenvalues $\gamma$ in the upper
(lower) half-plane. Because $A$ and $K$ are real-valued, we have $N_{c^+} = N_{c^-}$.
\end{itemize}

We are now ready to reformulate Theorem 1 from \cite{ChP10}.

\begin{Theorem}\cite{ChP10}
\label{equality 1} Under assumptions (H1)--(H4), for a small positive number $\delta$,
eigenvalues of the generalized eigenvalue problem
(\ref{shifted_eigenvalue_problem}) are counted as follows:
\begin{eqnarray}
\label{negative-index-A} N_p^- + N_n^0 + N_n^+ + N_{c^+} & = &
{\rm dim}({\cal H}_{A_{\delta}}^-), \\ \label{negative-index-K}
N_n^- + N_n^0 + N_n^+ + N_{c^+} & = & {\rm dim}({\cal H}_K^-).
\end{eqnarray}
\end{Theorem}

To apply Theorem \ref{equality 1} to the count of isolated
and embedded eigenvalues in the stability problem (\ref{spectrum}),
we recall from (\ref{count-K}) and (\ref{count-A-delta}) that
${\rm dim}({\cal H}_K^-) = n(\mathcal{L})$ and
${\rm dim}({\cal H}_{A_{\delta}}^-) = n(\mathcal{L})$. At the same time,
definition of $N_n^0$ yields $N_n^0 = n_0$, where $n_0$ is introduced
in Lemma \ref{lemma-negative-A}. Using these counts, we rewrite equalities (\ref{negative-index-A})
and (\ref{negative-index-K}) in the more explicit form:
\begin{eqnarray}
\label{negative-index-A-new} N_p^- + N_n^+ + N_{c^+} & = &
n(\mathcal{L})-n_0, \\ \label{negative-index-K-new}
N_n^- + N_n^+ + N_{c^+} & = & n(\mathcal{L}) - n_0.
\end{eqnarray}

We now need to compute numbers $N_p^-$, $N_n^-$, $N_n^+$, and $N_{c^+}$
for real and complex eigenvalues of the generalized eigenvalue problem (\ref{generalized_eigenvalue_problem}),
which are related to real, imaginary, and complex eigenvalues of the spectral problem (\ref{spectrum}).
Note that the imaginary eigenvalues of the spectral problem (\ref{spectrum})
may be embedded into the continuous spectrum of the operator $\partial_x \mathcal{L}$.

We recall here again the Hamiltonian symmetry, that is, if $\lambda \neq 0$ is a simple eigenvalue
of the spectral problem (\ref{spectrum}), then $-\lambda$ is also a simple eigenvalue
of the spectral problem (\ref{spectrum}) and both $\lambda$ and $-\lambda$ correspond to
the same double eigenvalue $\gamma = -\lambda^2 \neq 0$ of
the generalized eigenvalue problem (\ref{generalized_eigenvalue_problem}).

\begin{Lemma}
\label{lemma-inversion} Let $\lambda_j \in \R_+$ and
$\tilde{\lambda}_j = -\lambda_j \in \R_-$ be simple eigenvalues of the
spectral problem (\ref{spectrum}) associated with the real-valued eigenvectors
$v_j$ and $\tilde{v}_j$ in $D(\partial_x \mathcal{L}) \cap \dot{H}^{-1}(\mathbb{R})$. Then, we have
\begin{eqnarray}
\label{orthogonality2} \langle \mathcal{L} v_j^{\pm}, v_j^{\pm} \rangle = \pm
2 \langle \mathcal{L} \tilde{v}_j, v_j \rangle, \quad \langle \mathcal{L} v_j^{\pm}, v_j^{\mp} \rangle = 0,
\end{eqnarray}
where $v_j^{\pm} = v_j \pm \tilde{v}_j$ are linearly independent.
\end{Lemma}

\begin{Proof}
We recall that the eigenvectors $v_j$ and $\tilde{v}_j$ for distinct simple eigenvalues
$\lambda_j$ and $\tilde{\lambda}_j$ are linear
independent, hence the linear combinations $v_j^+$ and $v_j^-$ are linearly independent.
Since $\lambda_j \neq 0$ and $v_j$ is real-valued,
we have
$$
\langle \mathcal{L} v_j, v_j \rangle = \frac{1}{\lambda_j} \langle \mathcal{L} v_j, \partial_x \mathcal{L} v_j \rangle = 0.
$$
Similarly, $\langle \mathcal{L} \tilde{v}_j, \tilde{v}_j \rangle = 0$.
The orthogonality relations (\ref{orthogonality2}) hold by
direct computations.
\end{Proof}

\begin{Lemma}
\label{lemma-orthogonality3} Let $\lambda_j \in i \R_+$ and
$\bar{\lambda}_j = -\lambda_j \in i\R_-$ be simple eigenvalues of the
spectral problem (\ref{spectrum}) associated with the eigenvectors
$v_j$ and $\bar{v}_j$ in $D(\partial_x \mathcal{L}) \cap \dot{H}^{-1}(\mathbb{R})$. Then, we have
\begin{eqnarray}
\label{orthogonality5} \langle \mathcal{L} v_j^{\pm}, v_j^{\pm} \rangle =
2 \langle \mathcal{L} v_j, v_j \rangle, \quad \langle \mathcal{L} v_j^{\pm}, v_j^{\mp} \rangle = 0,
\end{eqnarray}
where $v_j^{\pm} = v_j \pm \bar{v}_j$ are linearly independent
and $\langle \mathcal{L} v_j, v_j \rangle$ is real.
\end{Lemma}

\begin{Proof}
Since operator $\mathcal{L}$ is real-valued, the eigenvector $v_j$ of
the spectral problem (\ref{spectrum}) with ${\rm Im}(\lambda_j) \neq 0$ has both
real and imaginary parts. Since $\lambda_j \neq 0$, we have
$$
\langle \mathcal{L} \bar{v}_j, v_j \rangle = \frac{1}{\lambda_j} \langle \mathcal{L} \bar{v}_j, \partial_x \mathcal{L} v_j \rangle = 0.
$$
Furthermore, since $\mathcal{L}$ is self-adjoint, we have
$\langle \mathcal{L} v_j, v_j \rangle = \langle \mathcal{L} \bar{v}_j, \bar{v}_j \rangle$.
The orthogonality equations (\ref{orthogonality5}) hold by direct computations.
\end{Proof}

\begin{Proof1}{\bf of Theorem \ref{main-theorem}.}
By symmetries of the linearized Hamiltonian system, each
eigenvalue $\gamma_j = - \lambda_j^2$
of the generalized eigenvalue problem
(\ref{generalized_eigenvalue_problem}) has a double multiplicity
compared to the eigenvalue $\lambda_j$ of the spectral problem (\ref{spectrum}).
From two linearly independent eigenvectors
$v_j^{\pm} \in D(\partial_x \mathcal{L}) \cap \dot{H}^{-1}(\R) \subset L^2(\R)$ constructed
in Lemmas \ref{lemma-inversion} and \ref{lemma-orthogonality3}, we obtain
two linearly independent eigenvectors $w_j^{\pm} = \partial_x^{-1} v_j^{\pm}
\in \mathcal{H}$ of the generalized eigenvalue problem (\ref{generalized_eigenvalue_problem})
\footnote{The same count holds in the case of complex eigenvalues $\lambda_j$ with
${\rm Re}(\lambda_j) \neq 0$ and ${\rm Im}(\lambda_j) \neq 0$. If $v_j$ and $\tilde{v}_j$
denote linearly independent eigenvectors of the spectral problem (\ref{spectrum}) for
complex eigenvalues $\lambda_j$ and $-\lambda_j$, then we can define two
linearly independent eigenvectors $w_j = \partial_x^{-1} v_j$ and $\tilde{w}_j = \partial_x^{-1} \tilde{v}_j$
of the generalized eigenvalue problem (\ref{generalized_eigenvalue_problem}) in $\mathcal{H}$
for the double eigenvalue $\gamma_j = -\lambda_j^2$. In the case of complex eigenvalues,
we do not care about the values of the quadratic form associated with the operator $\mathcal{L}$ computed at the eigenvectors.}.

By the orthogonality condition (\ref{orthogonality2}), we have $N_n^- = N_p^-$ for
a negative eigenvalue $\gamma_j = -\lambda_j^2$ corresponding to two real eigenvalues
$\lambda_j$ and $-\lambda_j$.  Since $N_n^- + N_p^- = 2 N_r$ because of the double multiplicity
of eigenvalues $\gamma_j$ compared to the multiplicity of eigenvalues $\lambda_j$,
we obtain $N_n^- = N_p^- = N_r$. Similarly, for a complex eigenvalue $\gamma_j = -\lambda_j^2$
corresponding to two complex eigenvalues $\lambda_j$ and $-\lambda_j$, we
count $N_{c^+} = 2 N_c$.

By the orthogonality condition (\ref{orthogonality5}), the double
multiplicity of the positive eigenvalue $\gamma_j = -\lambda_j^2$ corresponding to
two imaginary eigenvalues $\lambda_j$ and $\bar{\lambda}_j = -\lambda_j$,
and the definition of $N_i^-$, we obtain $N_n^+ = 2 N_i^-$.
The count (\ref{closure-relation1}) follows equivalently from either equality
(\ref{negative-index-A-new}) or (\ref{negative-index-K-new})\footnote{From comparison between
(\ref{negative-index-A}) and (\ref{negative-index-K}), which is justified by Lemma
\ref{lemma-ess-spectrum} and the equality $N_p^- = N_n^-$,
we obtain ${\rm dim}(\mathcal{H}_{A_{\delta}}^-) = {\rm dim}(\mathcal{H}_{K}^-) = n(\mathcal{L})$,
which yields the second independent proof of Lemma \ref{lemma-A-delta-K}.}.
\end{Proof1}

\begin{Remark}
The count of eigenvalues provided by the equality (\ref{negative-index-K}) in Theorem \ref{equality 1}
is a sufficient tool to prove Theorem \ref{main-theorem} since it follows from definitions that
${\rm dim}({\cal H}_K^-) = n(\mathcal{L})$, $N_n^0 = n_0$, whereas it follows from
Lemmas \ref{lemma-inversion} and \ref{lemma-orthogonality3} that $N_n^- = N_r$,
$N_n^+ = 2 N_i^-$, and $N_{c^+} = 2 N_c$. Along this avenue, the count
of negative eigenvalues of operators $\mathcal{M}$, $A$, and $A_{\delta}$ in Lemmas
\ref{lemma-operator-M}, \ref{lemma-negative-A}, and \ref{lemma-A-delta-K} (which is not
so easy to prove) is redundant and unnecessary.\label{remark}
\end{Remark}

\section{Generalization of Theorem \ref{main-theorem} for a periodic nonlinear wave}

We shall now take Remark \ref{remark} into account for an easy proof of the instability index count
for periodic waves in the KdV-type evolution equations. These instability index counts
were reported in \cite{BronKap,DecKap,KapHar} by means of much longer and different analysis.

We now consider a $2 L$-periodic nonlinear wave $\phi$ in a subspace $\mathcal{X}$ of
Hilbert space $L^2_{\rm per}(-L,L)$ equipped with the inner product
$\langle \cdot, \cdot \rangle$\footnote{Note that we do not change notations for
the inner product compared to the case $L^2(\R)$ but understand that the integration
is now performed on $[-L,L]$.} and an induced norm $\| \cdot \|$.
The spectral stability of $\phi$ is still determined by the
spectral problem (\ref{spectrum}), where $\mathcal{L} := E''(\phi)$ is a self-adjoint
real-valued operator with a dense domain $D(\mathcal{L})$ in $L^2_{\rm per}(-L,L)$.
We assume that $\mathcal{L}$ has a compact resolvent, so that the
spectrum of $\mathcal{L}$ in $L^2_{\rm per}(-L,L)$ is purely discrete.
We reinforce assumptions (H2) and (H3) in the slightly modified form:

\begin{itemize}
\item[({\bf H1})] The spectrum of $\mathcal{L}$ is purely discrete and
includes a finite number $n(\mathcal{L})$ of negative eigenvalues with
eigenvectors in $D(\mathcal{L}) \subset L^2_{\rm per}(-L,L)$.

\item[({\bf H2})] ${\rm Ker}(\mathcal{L}) = {\rm span}\{f_0\}$ with
$f_0 \in D(\mathcal{L}) \cap D(\partial_x \mathcal{L} \partial_x) \cap \dot{H}_{\rm per}^{-1}(-L,L)$,
so that $\phi_0 = \partial_x^{-1} f_0 \in L^2_{\rm per}(-L,L)$.
\end{itemize}

In addition, we note that ${\rm Ker}(\partial_x) = {\rm span}\{1\} \subset L^2_{\rm per}(-L,L)$ and define
the matrix $\mathcal{D}$ as follows:
\begin{equation}
\label{matrix-D}
\mathcal{D} = \left[ \begin{array}{cc} \langle \mathcal{L}^{-1} \phi_0,\phi_0 \rangle & \langle \mathcal{L}^{-1} \phi_0, 1 \rangle \\
\langle \mathcal{L}^{-1} \phi_0, 1 \rangle & \langle \mathcal{L}^{-1} 1, 1 \rangle \end{array} \right].
\end{equation}
Note that matrix $\mathcal{D}$ has finite elements because
${\rm span}\{1,\phi_0\} \perp {\rm Ker}(\mathcal{L})$ as
it follows from the orthogonality conditions
$\langle f_0, \phi_0 \rangle = \langle \partial_x \phi_0, \phi_0 \rangle = 0$
and $\langle f_0, 1 \rangle = \langle \partial_x \phi_0, 1 \rangle = 0$.
We modify now assumption (H4) as follows:
\begin{itemize}
\item[({\bf H3})] Matrix $\mathcal{D}$ is invertible.
\end{itemize}

With the previous definitions of $N_r$, $N_c$, and $N_i^-$,
the following theorem gives a modification of Theorem \ref{main-theorem}
for a periodic nonlinear wave.

\begin{Theorem}
\label{theorem-periodic-wave}
Assume ({\bf H1})--({\bf H3}). Then,
\begin{equation}
\label{closure-relation2} N_r + 2 N_c + 2 N_i^- = n(\mathcal{L}) - n(\mathcal{D}),
\end{equation}
where $n(\mathcal{D})$ is the number of negative eigenvalues of the matrix $\mathcal{D}$.
\end{Theorem}

\begin{Proof}
We extend the spectral problem (\ref{spectrum}) to the system of two coupled equations:
\begin{equation}
\label{coupled-eigenvalue}
\mathcal{M} w = -\lambda v, \quad \mathcal{L} v = \lambda w,
\end{equation}
where $\mathcal{M} = -\partial_x \mathcal{L} \partial_x$, $v \in D(\partial_x \mathcal{L}) \cap \dot{H}_{\rm per}^{-1}(-L,L)
\subset L^2_{\rm per}(-L,L)$,
and $w \in D(\partial_x \mathcal{L} \partial_x) \subset L^2_{\rm per}(-L,L)$. The equivalence
of simple eigenvalues of the coupled system (\ref{coupled-eigenvalue}) and those of the spectral
problem (\ref{spectrum}) is proved similarly to Proposition \ref{equivalence-lemma}.

By assumptions ({\bf H1}) and ({\bf H2}), the zero eigenvalue of $\mathcal{L}$ is isolated and simple.
Let $P$ be the orthogonal projection from $L^2_{\rm per}(-L,L)$ to $[{\rm span}\{f_0\}]^{\perp} \subset L^2_{\rm per}(-L,L)$.
By a procedure that is similar to (\ref{equation-w}), (\ref{equation-u}), and (\ref{equation-v}),
we obtain the generalized eigenvalue problem for a nonzero eigenvalue $\gamma \neq 0$:
\begin{equation}
\label{generalized-eigenvalue}
A w = \gamma K w, \quad w \in \mathcal{H},
\end{equation}
where $A := P \mathcal{M} P$, $K := P \mathcal{L}^{-1} P$, $\gamma := - \lambda^2$, and
$\mathcal{H} := D(\mathcal{M}) \cap [{\rm span}\{f_0\}]^{\perp} \subset L^2_{\rm per}(-L,L)$.
The equivalence of double eigenvalues of the generalized eigenvalue problem (\ref{generalized-eigenvalue})
and pairs of simple eigenvalues of the coupled system (\ref{coupled-eigenvalue}) is proved
similarly to Proposition \ref{proposition-generalized-eigenvalue}.

The spectrum of $\mathcal{M}$ is purely discrete but the
zero eigenvalue of $\mathcal{M}$ is now double since ${\rm Ker}(\mathcal{M}) = {\rm span}\{1,\phi_0\}$.
Therefore, for a small positive number $\delta$, assumptions of Theorem \ref{equality 1}
are satisfied and the equality (\ref{negative-index-K}) takes the form:
\begin{equation}
\label{count-1}
N_n^- + N_n^0 + N_n^+ + N_{c^+} = {\rm dim}(\mathcal{H}_K^-).
\end{equation}
By construction of $K = P \mathcal{L}^{-1} P$ and
$\mathcal{H} = D(\mathcal{M}) \cap [{\rm span}\{f_0\}]^{\perp} \subset L^2_{\rm per}(-L,L)$,
we have ${\rm dim}(\mathcal{H}_K^-) = n(\mathcal{L})$. On the other hand,
$N_n^0$ denotes algebraic multiplicity of zero eigenvalues of the generalized eigenvalue problem
(\ref{generalized-eigenvalue}) whose generalized eigenvectors are associated to non-positive values of the
quadratic form $\langle K \cdot,\cdot \rangle$. Since ${\rm Ker}(\mathcal{M}) = {\rm span}\{1,\phi_0\}$
and the matrix $\mathcal{D}$ has no zero eigenvalue by assumption ({\bf H3}),
we have $N_n^0 = n(\mathcal{D})$.

From analysis identical to Lemmas \ref{lemma-inversion} and \ref{lemma-orthogonality3},
we also obtain $N_n^- = N_r$, $N_n^+ = 2 N_i^-$, and $N_{c^+} = 2N_c$,
hence the count (\ref{count-1}) yields the instability index count (\ref{closure-relation2})
and the theorem is proven.
\end{Proof}

\begin{Remark}
Equality (\ref{negative-index-A}) in Theorem \ref{equality 1} can also be used for the correct
instability index count but this task would require the count of negative eigenvalues of operators
$\mathcal{M}$, $A$, and $A_{\delta}$ similar to that in Lemmas
\ref{lemma-operator-M}, \ref{lemma-negative-A}, and \ref{lemma-A-delta-K},
which would result in the formula ${\rm dim}(\mathcal{H}_{A + \delta K}^-) = n(\mathcal{L})$ for a small
positive number $\delta$\footnote{The identity ${\rm dim}(\mathcal{H}_{A + \delta K}^-) = n(\mathcal{L})$
follows from the identity $n(A) = n(P \mathcal{M} P)
= n(\mathcal{L}) - n(\mathcal{D})$, which should hold despite the fact that the projection operator $P$ is
defined by the orthogonal complement of the one-dimensional
subspace ${\rm Ker}(L) = {\rm span}\{f_0\}$. The corresponding argument goes as follows. To study $n(P\mathcal{M}P)$,
we introduce a Lagrange multiplier $\nu$ and set up the self-adjoint spectral problem $\mathcal{M} w = \mu w + \nu f_0$
with the orthogonality condition
$\langle f_0, w \rangle = 0$. Since $\mathcal{M} = -\partial_x \mathcal{L} \partial_x$, $f_0 = \partial_x \phi_0$,
and ${\rm Ker}(\partial_x) = {\rm span}\{1\}$, we set $w = \partial_x g$, integrate in $x$, and
obtain the non-self-adjoint spectral problem $\mathcal{L} (-\partial_x^2) g = \mu g + \nu \phi_0 + \chi$
with two Lagrange multipliers $\nu$ and $\chi$, under the constraints $\langle \phi_0, \partial_x^2 g \rangle = 0$
and $\langle 1, \partial_x^2 g \rangle = 0$. Smoothing it with a positive parameter $\epsilon$
and setting $g = (\epsilon - \partial_x^2)^{-1} h$, we end up with the self-adjoint spectral problem
$\mathcal{L} h = \mu (\epsilon - \partial_x^2)^{-1} h + \nu \phi_0 + \chi$ under the constraints
$\langle \phi_0, h \rangle = 0$ and $\langle 1, h \rangle = 0$, which can be studied
with the standard technique (Theorem 4.1 in \cite{Pel-book}).}. For the spectral problem associated with the KdV-type evolution
equation, this equality is redundant because the spectral problem for the coupled system (\ref{coupled-eigenvalue})
is a squared version of the original spectral problem (\ref{spectrum}).
\end{Remark}

\begin{Remark}
Eigenvectors of the spectral problem (\ref{spectrum}) in $L^2_{\rm per}(-L,L)$ for a nonzero eigenvalue
$\lambda$ are orthogonal to ${\rm Ker}(\mathcal{M}) = {\rm span}\{1,\phi_0\}$. Therefore, we can introduce a constrained space
$$
\tilde{\mathcal{H}} := D(\mathcal{L}) \cap [{\rm span}\{1,\phi_0\}]^{\perp} \subset L^2_{\rm per}(-L,L)
$$
and a projection operator $\tilde{P} : L^2_{\rm per}(-L,L) \to \tilde{\mathcal{H}}$ to
reformulate the spectral problem (\ref{spectrum}) as a linearized Hamiltonian
system with an invertible symplectic matrix:
\begin{equation}
\label{operator-R-S}
\mathcal{L}_p v = \lambda J_p v, \quad v \in \tilde{\mathcal{H}},
\end{equation}
where $\mathcal{L}_p := \tilde{P} \mathcal{L} \tilde{P}$ and
$J_p := \tilde{P} \partial_x^{-1} \tilde{P}$.
Standard analysis in constrained Hilbert spaces (Theorem 4.1 in \cite{Pel-book})
shows that $n(\mathcal{L}_p) = n(\mathcal{L}) - n(\mathcal{D})$. Applying now
the general instability index count in the linearized Hamiltonian systems
with an invertible symplectic operator \cite{KKS04}, one can immediately obtain
the instability index count formula (\ref{closure-relation2}). This proof of the instability
index count for the nonlinear periodic waves in Hamiltonian systems was introduced
by Haragus \& Kapitula \cite{KapHar}.
\end{Remark}

Let us show how to recover the correct count of eigenvalues for the example of the focusing
modified KdV equation
\begin{equation}
\label{mKdV}
u_t + 3 u^2 u_x + u_{xxx} = 0.
\end{equation}
Travelling periodic waves in the form
$u = \phi(x-ct)$ satisfies the differential equation
\begin{equation}
\label{stationary-mKdV}
\phi'' = c \phi - \phi^3,
\end{equation}
where the constant of integration is chosen to be zero.
Two families of nonlinear periodic waves were considered by
Deconinck \& Kapitula \cite{DecKap} in the explicit form:
\begin{eqnarray}
\label{solution-1}
\phi(x) & = & \sqrt{2} k {\rm dn}(x,k), \quad c = 2 - k^2,\\
\label{solution-2}
\phi(x) & = &  \sqrt{2} k {\rm cn}(x,k), \quad c = -1 + 2k^2,
\end{eqnarray}
where ${\rm dn}$ and ${\rm cn}$ are Jacobi's elliptic functions and
the period $L = 4 K(k)$ is given by the complete elliptic integral of the first kind
for a fixed $k \in (0,1)$.

Since $\mathcal{L} := -\partial_x^2 + c - 3 \phi^2(x)$ with
$\mathcal{L} \partial_x \phi = 0$ and $\mathcal{L} \partial_c \phi = -\phi$,
assumption ({\bf H2}) is satisfied with $\phi_0 = \phi$. In addition,
by scaling invariance of the stationary modified KdV equation (\ref{stationary-mKdV}), one can
obtain that
\begin{eqnarray}
\langle \mathcal{L}^{-1} \phi,\phi \rangle = -\frac{1}{2} \frac{d}{dc} \| \phi \|^2 < 0
\end{eqnarray}
and
\begin{eqnarray}
\langle \mathcal{L}^{-1} \phi,1 \rangle = -\frac{d}{dc} \int_{-L}^L \phi(x) dx = 0.
\end{eqnarray}

Let us denote
\begin{equation}
F(k) := \langle \mathcal{L}^{-1} 1, 1 \rangle = \frac{1}{2L} \int_{-L}^L \mathcal{L}^{-1}(1) dx.
\end{equation}
For the ${\rm dn}$-wave, explicit computations in \cite{DecKap} show that
$n(\mathcal{L}) = 1$ and $F(k) > 0$ for all $k \in (0,1)$. Therefore,
the instability index count (\ref{closure-relation2}) with $n(\mathcal{L}) = n(\mathcal{D}) = 1$ yields
spectral stability of the ${\rm dn}$-wave for all $k \in (0,1)$.

For the ${\rm cn}$-wave, explicit computations in \cite{DecKap} show that
$n(\mathcal{L}) = 2$, whereas there is $k_* \approx 0.909$ such that $F(k) < 0$ for $0 < k < k^*$
and $F(k) > 0$ for $k^* < k < 1$. Therefore, the instability index count (\ref{closure-relation2})
yields spectral stability of the ${\rm cn}$-wave for $k \in (0,k_*)$
with $n(\mathcal{L}) = n(\mathcal{D}) = 2$ and spectral instability for $k \in (k_*,1)$
with $n(\mathcal{L}) = 2$, $n(\mathcal{D}) = 1$, and $N_r = 1$.

\section{Conclusion}

Although the instability index count is now well established both
for solitary waves and periodic waves of the KdV-type evolution equations,
there are still many directions of further development in the stability theory of nonlinear waves.
In particular, Boussinesq equations involve the spectral problem
(\ref{spectrum}) associated with a matrix differential operator $\mathcal{L}$ \cite{Lin}.
Although this matrix operator can be mapped to a self-adjoint form
with a similarity transformation, it becomes difficult to transform both
operators $\mathcal{L}$ and $\mathcal{M}$ in the coupled system (\ref{coupled-problem})
to the self-adjoint form. More direct approaches to the stability of
solitary waves in Boussinesq systems can be found in recent works \cite{Stef1,Stef2}.
The same complication may also occur in the system of coupled KdV-type equations.

Further extensions of the instability index count involves quadratic operator
pencils, as well as general polynomial pencils, where the count of unstable eigenvalues
become less precise. Works in this direction can be found in \cite{Bronski,ChPel,Kollar}.

\end{document}